\documentclass[landscape, oneside, 12pt]{amsart}
\usepackage{xfrac}
\usepackage[margin=7cc]{geometry}
\usepackage{mathscinet}
\usepackage{amssymb}
\usepackage{tikz}
\usepackage{doi}
\usepackage{placeins}
\usepackage[doublespacing]{setspace}

\DeclareMathOperator{\Dense}{Dense}
\DeclareMathOperator{\meagre}{\mathcal{M}}
\DeclareMathOperator{\lmz}{\mathcal{N}}
\DeclareMathOperator{\smz}{\mathcal{SN}}
\DeclareMathOperator{\add}{add}
\DeclareMathOperator{\cov}{cov}
\DeclareMathOperator{\non}{non}
\DeclareMathOperator{\cof}{cof}

\newcommand{\Pot}{\mathcal{P}}
\newcommand{\swhalf}{\mathfrak{s}^w_{\sfrac{1}{2}}}
\newcommand{\QQ}{\mathbb{Q}}
\newcommand{\ZZ}{\mathbb{Z}}
\newcommand{\RR}{\mathbb{R}}
\newcommand{\cf}{cf.}
\newcommand{\subsetmf}{\subset^*}
\newcommand{\subsetmnwd}{\subset_\text{\upshape nwd}}
\newcommand{\leqmf}{\leqslant^*}
\newcommand{\dfeq}{:=}
\newcommand{\wlg}{without loss of generality}
\newcommand{\ifff}{if and only if}

\newcommand{\headline}[1]{\text{\bfseries\upshape #1}}

\newcommand{\seq}[2]{\langle #1 \, \mid \, #2\rangle}

\newcommand{\eqfr}{\mathfrak{e\mspace{-0.5mu}
q}}
\newcommand{\spair}{\mathfrak{s}_\text{\itshape pair}}
\newcommand{\shalfinf}{\mathfrak{s}_{\sfrac{1}{2}}^\infty}
\newcommand{\shalfpme}{\mathfrak{s}_{\sfrac{1}{2}\pm\varepsilon}}
\newcommand{\sss}{\text{\upshape\bfseries\ss}}
\newcommand{\eubd}{\mathfrak{e}_\text{ubd}}

\newcommand{\absv}[1]{\vert #1\vert}
\newcommand{\card}[1]{\lvert #1\rvert}

\newcommand{\soast}[2]{\left\{#1\mid #2\right\}}
\newcommand{\bsoast}[2]{\big\{ #1 \, \big\vert \, #2\big\}}

\newcommand{\bgsoast}[2]{\bigg\{ #1 \, \bigg\vert \, #2\bigg\}}

\newcommand{\soafft}[2]{{}^{#1}  #2}
\newcommand{\opair}[2]{{( #1, #2)}}

\DeclareMathSizes{10}{10}{8}{7}
\DeclareMathSizes{11}{11}{9}{8}
\DeclareMathSizes{12}{12}{10}{9}

\theoremstyle{thmstyleone}
\newtheorem{defi}{Definition}[subsection]
\newtheorem{theo}[defi]{Theorem}
\newtheorem{lemm}[defi]{Lemma}
\newtheorem{coro}[defi]{Corollary}
\newtheorem{claim*}{Claim}

\newenvironment{proofofclaim}{{\itshape Proof of Claim.}}{\hfill\itshape(End of Proof of Claim.)}

\begin{document}
\title[Comparing Cardinal Characteristics]{Calligraphy Concerning \\
Casually Compiled Cardinal Characteristics Comparisons}

\author{Thilo Weinert}
\email{thilo.weinert@univie.ac.at}


\begin{abstract}
The paper establishes several inequalities between cardinal characteristics of the continuum. In particular, it is shown that the partition splitting number is not larger than the uniformity of the meagre ideal; not all sets of reals having the cardinality of an the $\varepsilon$-almost bisecting number are of strong measure zero; no fewer sets of strong measure zero than indicated by the statistically reaping number suffice to cover the reals; the pair-splitting number is not smaller than the evasion number; and the subseries number is neither smaller than the pair-splitting number nor than the minimum of the unbounding number and the unbisecting number. Moreover, a diagram putting these results into context is provided and a brief historical account is given.\end{abstract}


\maketitle

\section{Preliminaries}

\subsection{Formalities}

\subsubsection{Acknowledgements}
The author thanks Andreas Blass, Shimon Garti, Michael Hru\v{s}\'{a}k, Ferdinand Ihringer, and Diana Montoya for enlightening correspondence by electronic mail and gratefully acknowledges support through the Singaporean Ministry of Education Tier 2 grant No. MOE2017-T2-2-125, the Fonds zur F\"orderung der wissenschaftlichen Forschung, Lise Meitner grant M 3037-N, as well as the Research Project of National relevance ``PRIN2022\_DIMONTE - Models, sets and classifications\footnote{realizzato con il contributo del progetto PRIN 2022 - D.D. n. 104 del 02/02/2022 – PRIN2022\_DIMONTE - Models, sets and classifications - Codice 2022TECZJA\_003 - CUP N. G53D23001890006. “Finanziato dall'Unione Europea – Next-GenerationEU – M4 C2 I1.1}
. For the purpose of open access, the author has applied a CC BY public copyright licence to any Author Accepted Manuscript version arising from this submission.

\subsubsection{Declarations}

The author has no relevant financial or non-financial interests to disclose.

The author has no conflicts of interest to declare that are relevant to the content of this article.

The author certifies that he have no affiliation with or involvement in any organisation or entity with any financial interest or non-financial interest in the subject matter or materials discussed in this manuscript.

The author has no financial or proprietary interest in any material discussed in this article.

\clearpage
\subsection{Introduction}

\subsubsection{Prologue}
Cardinal characteristics of the continuum are infinite cardinals, usually uncountable and no larger than the continuum, given by a specific definition. They are often derived from considerations not belonging to set theory proper and sometimes from investigations in combinatorial set theory for their own sake. There are two kinds of possible results about them: Often one can show metamathematically that it is consistent with the axioms of set theory to assume that one particular characteristic is smaller than another; sometimes one can straightforwardly show that one is no larger than another. This paper contains six results of the latter sort, in general related to work of Blass, Hru\v{s}\'{a}k, Laflamme, Meza-Alcántara, A. Miller, Minami, Raghavan and Stepr\'-{a}ns, \cite{989BL0, 994B0, 004BHH0, 006MS0, 010M0, 010HMM0, 020R0} and is particularly inspired by the more recent work of J. Brendle, Brian, Halbeisen, Hamkins, Klausner, Lischka, and Shelah, \cite{023BHKLS0, 019BBH0}. Its purpose is twofold. On the one hand it shall disseminate the results presented among those specialising---among possibly other fields---in the theory of cardinal characteristics of the continuum. On the other hand, it shall give a broader audience a small overview over parts of the subject and serve as a possible starting point for enquiries. While thus having aspects of a survey, it does not aim to be an exposition. For a more expository treatment, see \cite{984D0, 995BJ0, 001B0, 010B0, 012H1}. For a thorough historical account, see \cite{012S1}.
Combining these objectives seems appropriate as, firstly, none of the results presented are of a metamathematical nature, so the reader does not need to know any mathematical logic, and, secondly, they are somewhat scattered 
within the subject of cardinal characteristics of the continuum. A Hasse diagram is provided on page \pageref{diagram}. It includes all cardinal characteristics involved in new results presented here---see the Summary on page \pageref{summary}---and then some. The general aim was to include as many cardinal characteristics as possible without rendering the resulting diagram incomprehensible. In particular, if a characteristic appears in the diagram, known relations to other characteristic in the diagram appear as well; moreover, lines are not crossing. These conditions explain why $\mathfrak{b}$ does appear in the diagram within both a minimum and a maximum but not on its own, as $\mathfrak{h}_\QQ \leqslant \mathfrak{b}$. No reader should infer that a cardinal characteristic is generally wider known or in some objective sense more interesting than another only because in contrast to the latter, the former appears in the diagram. A line between two characteristics means an inequality applies with the characteristic higher up on the page being no smaller than the one further down on it.

All characteristics mentioned in this paper have been introduced by someone else elsewhere. Whenever a characteristic is defined, a reference to a survey or one or several papers concerning the characteristic is given.

All proofs presented in this paper are by contradiction and written unwaveringly in this manner. This is to say, that the contradiction is not reached before the end of the proof. Their style might strike the readers as terse but the author hopes that they will (grow to) like it. It may be less suited to give readers a rough idea but it probably is all the more suited to thoroughly check for correctness which, in the absence of a coauthor, is an important feature.

\subsubsection{History}
It only gradually became a custom to name cardinal characteristics when they appeared in a theorem. This is, by the way, the reason that the legend on page \pageref{legend} contains two columns with references: One to a place in which an inequality between characteristics can be said to have appeared earliest in a history of ideas, and one to a source where explicit definitions and proofs in modern terminology can be found, would this not be possible in the the original source. Many characteristics were named by van Douwen in \cite{984D0} but appeared in mathematics earlier. It later became customary to depict characteristics thus named in a Hasse diagram called \emph{van Douwen's diagram}. 

The unbounding number $\mathfrak{b}$, \cf\ Definition \ref{unbounding number}, goes back to Rothberger and Sierpi\'nski, \cite{939R0,939S0}, the splitting number $\mathfrak{s}$, \cf\ Definition \ref{splitting number}, to Booth, \cite{974B3} and the distributivity number $\mathfrak{h}$, \cf\ Definition \ref{distributivity number}, to Balcar, Pelant, and Simon, \cite{980BPS0}. The study of cardinal characteristics related to measure and category became oriented around what, following Fremlin, \cite{984F0}, became known as Cichoń's diagram. Again, investigations had begun before that, at the latest with A. Miller's \cite{981M0}. The groupwise density number $\mathfrak{g}$, \cf\ Definition \ref{groupwise density number} was introduced by Blass and Laflamme, \cite{989BL0}. In \cite{994B0}, Blass defined the \emph{evasion number} $\mathfrak{e}$, \cf\ Definition \ref{evasion number} and the linear evasion number $\mathfrak{e}_\ell$, and used them to analyse what he termed the \emph{Specker phenomenon} thus continuing earlier investigations by Eda, \cite{983E0}.
A subgroup $G$ of the Baer-Specker-group $\ZZ^\omega$,  exhibits the \emph{Specker phenomenon}
if it contains a sequence $\seq{a_i}{i < \omega}$ of linearly independent elements such that for every homomorphism
$h : G \longrightarrow \ZZ$ the set
$\soast{i < \omega}{h(a_i) \ne 0}$ is finite.
Blass showed that there is a group no larger than $\min(\mathfrak{b}, \mathfrak{e})$ exhibiting the Specker phenomenon and that the smallest group exhibiting it can be no smaller than $\mathfrak{e}_\ell$. Brendle, \cf\ \cite{996BS0}, went on to show that the smallest group exhibiting it is in fact of cardinality $\mathfrak{e}_\ell$ which he proved to equal $\min(\mathfrak{b}, \mathfrak{e})$.

Another development starting in the nineties of the last century with a paper by Majcher, \cite{990M1}, is the dualisation of van Douwen's diagram. The notion of duality here is the same as in the dual form of Ramsey's Theorem, proved by Carlson and Simpson, \cf\ \cite{984CS0}: One switches from subsets of a set to partitions of it and as subsets are induced by ranges of injections, partitions are induced by domains of surjections. These dual characteristics were subsequently analysed by several people, \cf e.g.\cite{000CKMW0}.
In \cite{010M0}, Minami analysed the dual splitting number and the dual reaping number $\mathfrak{r}_d$ and to that end introduced the notions of pair-splitting and pair-reaping.

The characteristic $\eqfr$ was introduced officially by A. Miller and Stepr\-{a}ns, \cite[p. 57, Def. 16]{006MS0}, but one could argue that it was in the air since Miller proved that it is a combinatorial description of the minimal size of a set of reals failing to have strong measure zero, \cf\ \cite[p. 97, Thm. 2.3]{981M0}. Note, however, that this characterisation holds for $\RR$ and the Cantor space but fails for the Baire space as discussed by Hru\v{s}\'{a}k, Wohofsky, and Zindulka, \cf\ \cite{016HWZ0}.

The cardinal characteristics in Definition \ref{starred characteristics} were first considered under different names in \cite{999BS0} and subsequently analysed in \cite{007HH0, 017RS0, 020R0}.

Another variation of van Douwen's diagram comes from studying the structure $\opair{\Dense(\QQ)}{\subsetmnwd}$ instead of $\opair{\Pot(\omega)}{\subsetmf}$. This was done by Balcar, Brendle, Cichoń, Hernández-Hernández, and Hru\v{s}\'{a}k, \cf\ \cite{001C1, 004BHH0, 006B2}.

Topological selection principles inspired the definition of the linear refinement number, \cite[Definition 61]{003T0} and the linear excluded middle number, \cite[Definition 1.2]{016MST0}.

Most recently, Blass, J. Brendle, Brian, Hamkins, Hardy, and P. Larson analysed the number of manipulations needed to change the limit of a conditionally convergent series, by respectively rearranging, \cf\ \cite{020BBBHHL0}, and dropping some summands, \cf\ \cite{019BBH0}. Meanwhile J. Brendle, Halbeisen, Klausner, Lischka, and Shelah investigated variants of the splitting and unsplitting numbers derived from requiring splitting to happen equitably in various ways, \cite{023BHKLS0}.

\subsection{Notation, Definitions, a Summary, and a Diagram}

\subsubsection{Notation and Definitions}

For a set $A$, we write $\Pot(A)$ for the power set of $A$; $[A]^\kappa$ for the families of subsets of $A$ having cardinality $\kappa$; and $[A]^{<\kappa}$ for the families of subsets of $A$ having fewer than $\kappa$ elements.
For sets $A$ and $B$, by writing $\soafft{A}{B}$, we refer to the collection of all functions from $A$ to $B$. We employ von-Neumann-ordinals, which means that every ordinal number is equal to the set of its predecessors. So $0$ denotes the number zero but also the empty set and $\omega$ denotes the smallest transfinite ordinal but also the set of all natural numbers.  Moreover, this means that for two ordinals $\alpha$ and $\beta$, the set $\alpha \setminus \beta$ is the half-open interval from $\beta$ (and if $\beta < \alpha$ including it) to $\alpha$ but excluding it. By writing $X \subset Y$, we mean that $X$ is a subset of $Y$, not necessarily a proper subset. Writing $X \subsetmf Y$ means that $Y \setminus X$ is finite while for sets $X$ and $Y$ of rational numbers, $X \subsetmnwd Y$ means that $Y \setminus X$ is nowhere dense in $\QQ$.

We denote by $\meagre$ the ideal of meagre (or first category) sets, by $\lmz$ the ideal of sets of Lebesgue measure $0$, and by $\smz$ the ideal of sets of strong measure zero. Recall that a set $X$ is of \emph{strong measure zero} if for every sequence $\seq{\varepsilon_i}{i < \omega}$ there is a sequence $\seq{I_n}{n < \omega}$ of open intervals with $I_n$ having length at most $\varepsilon_n$ for every natural number $n$ such that $\displaystyle X \subset \bigcup_{n < \omega} I_n$. Finally, $\mathcal{Z}_0$ denotes the ideal of sets of natural numbers that have asymptotic density $0$.

\begin{defi}[{\cite[p. 115]{984D0}, \cite[2.2 Def.]{010B0}, \cite[Section 8]{012H1}}]
\label{unbounding number}
The \emph{unbounding number} (sometimes called the \emph{bounding number}) $\mathfrak{b}$ is the minimal cardinality of a family $\mathcal{F}$ of functions $f \in \soafft{\omega}{\omega}$ such that there is no function $g \in \soafft{\omega}{\omega}$ satisfying $f \leqmf g$ for all $f \in \mathcal{F}$.
\end{defi}

A family $\mathcal{F} \subset [\omega]^\omega$ is called \emph{open} if for every $X \in \mathcal{F}$ we have $Y \in \mathcal{F}$ for all $Y \subsetmf X$. It is called \emph{dense} if for every $X \in [\omega]^\omega$ there is a $Y \in \Pot(X) \cap \mathcal{F}$.

A family $\mathcal{G} \subset [\omega]^\omega$ is \emph{groupwise dense} if it is open and for every interval partition $\seq{J_i}{i < \omega}$ there is an infinite set $A$ of natural numbers such that $\displaystyle\bigcup_{i \in A} J_i \in \mathcal{G}$.

\begin{defi}[{\cite[p. 53]{989BL0}, \cite[6.26 Def.]{010B0}}]
\label{groupwise density number}
The \emph{groupwise density number} $\mathfrak{g}$ is the minimal cardinality of a collection $\mathcal{C}$ of groupwise dense families such that $\bigcap\mathcal{C} = 0$.
\end{defi}

\begin{defi}[{\cite[6.5 Def.]{010B0}, \cite[Section 8]{012H1}}]
\label{distributivity number}
The \emph{distributivity number} $\mathfrak{h}$ is the minimal cardinality of a collection $\mathcal{C}$ of dense open families of sets of natural numbers such that $\bigcap\mathcal{C} = 0$.
\end{defi}

Following \cite{004BHH0}, we will denote the family of all dense sets of rational numbers by $\Dense(\QQ)$. A family $\mathcal{F} \subset \Dense(\QQ)$ is called $\QQ$-\emph{open} if for every $X \in \mathcal{F}$ we have $Y \in \mathcal{F}$ for all $Y \subsetmnwd X$. Moreover, a family $\mathcal{D} \subset \Dense(\QQ)$ is called $\QQ$-\emph{dense} if for every $X \in \Dense(\QQ)$ there is a $Y \in \Pot(X) \cap \mathcal{D}$.
 
\begin{defi}[{\cite[Section 3]{004BHH0}}]
The \emph{distributivity number} of $\opair{\Dense(\QQ)}{\subsetmnwd}$, written $\mathfrak{h}_\QQ$ is the minimal cardinality of a collection $\mathcal{C}$ of $\QQ$-dense $\QQ$-open families of sets of rationals such that $\bigcap\mathcal{C} = 0$.
\end{defi}

\begin{defi}
For infinite sets $S$ and $X$ of natural numbers and $\varepsilon \in ]0, \sfrac{1}{2}[$, we say
\begin{align*}
S \text{ \emph{splits} } X & \text{ if both } X \cap S \text{ and } X \setminus S \text{ are infinite,} \\
S \text{ \emph{bisects} } X & \text{ if } \lim_{n \nearrow \infty} \frac{\card{S \cap X \cap n}}{\card{X \cap n}} = \frac{1}{2}, \\
S \text{ \emph{infinitely often\footnotemark bisects} } X & \text{ if } \frac{\card{S \cap X \cap n}}{\card{X \cap n}} = \frac{1}{2} \text{ for infinitely many } n < \omega, \\
S\ \varepsilon\text{\emph{-almost bisects} } X & \text{ if } \bgsoast{n < \omega}{\frac{\card{S \cap X \cap n}}{\card{X \cap n}} \notin \left] \frac{1}{2} - \varepsilon, \frac{1}{2} + \varepsilon\right[} \text{ is finite.}
\end{align*}
\footnotetext{or cofinally}
\end{defi}

\begin{defi}
For sets $X, Y \in \Dense(\QQ)$ we say that $X$ $\QQ$\emph{-splits} $Y$ if both $Y \cap X$ and $Y \setminus X$ are dense in $\QQ$. A family $\mathcal{F}$ of dense sets of rational numbers is called $\QQ$\emph{-splitting} if for every dense set $X$ of rational numbers there is a member of $\mathcal{F}$ which $\QQ$-splits $X$.
\end{defi}

A set $X$ of natural numbers is called \emph{moderate} if
\begin{align*}
0 < \liminf_{n \nearrow \infty} \frac{\card{X \cap n}}{n} \leqslant \limsup_{n \nearrow \infty} \frac{\card{X \cap n}}{n} < 1.
\end{align*}
We also say that a set $A$ of pairs of natural numbers is \emph{unbounded} if
$A \not\subset [\omega]^2 \setminus [\omega \setminus k]^2$ for every natural number $k$. Moreover, an infinite set $X$ of natural numbers \emph{pair-splits} $A$ if $A \setminus ([X]^2 \cup [\omega \setminus X]^2)$ is infinite. 

\begin{defi}
A family $\mathcal{F}$ of infinite sets of natural numbers is called
\begin{align*}
\text{\emph{splitting} } & \text{if for every infinite set } X \text{ of natural numbers} \\
& \text{there is a member of } \mathcal{F} \text{ splitting } X, \\
\text{\emph{pair-splitting} } & \text{if for every unbounded } A \subset [\omega]^2 \\
& \text{there is a member of } \mathcal{F} \text{ pair-splitting it.} \\
\text{\emph{unbisected} } & \text{if there is no set } S \text{ of natural numbers} \\
& \text{bisecting every member of } \mathcal{F}, \\
\text{\emph{statistically refining\footnotemark} } & \text{if there is no moderate } S \in [\omega]^\omega \text{ such that} \\
&  \lim_{n \nearrow \infty} \frac{n\cdot\card{S \cap X \cap n}}{\card{S \cap n}\cdot\card{X \cap n}} = 1 \text{ for all } X \in \mathcal{F},\\
\text{\emph{bisecting} } & \text{if for every infinite set } X \text{ of natural numbers} \\
& \text{there is a member of } \mathcal{F} \text{ bisecting } X, \\
\text{\emph{cofinally bisecting} } & \text{if for every } X \in [\omega]^\omega \\
& \text{there is a member of } \mathcal{F} \text{ cofinally bisecting } X, \\
\varepsilon\text{\emph{-almost bisecting} for an } \varepsilon \in ]0, \sfrac{1}{2}[ & \text{ if for every } X \in [\omega]^\omega \\
& \text{there is a member of } \mathcal{F} \ \varepsilon\text{-almost bisecting } X.
\end{align*}
\footnotetext{or statistically reaping}
\end{defi}

The following notion of splitting number was introduced by Raghavan.

\begin{defi}[{\cite[Def. 3]{020R0}}]
Let $P$ be an infinite partition of $\omega$ and $a$ an infinite set of natural numbers. We say that $P$ \emph{splits} $a$ if $a \cap x$ is infinite for every $x \in P$. An family $\mathcal{F}$ of partitions of $\omega$ is called a \emph{splitting family of partitions} if for every $a \in [\omega]^\omega$ there is a $P \in \mathcal{F}$ splitting $a$. The \emph{partition splitting number} $\mathfrak{s}(\mathfrak{pr})$ is the minimal cardinality of a splitting family of partitions.
\end{defi}

It is currently open whether $\mathfrak{s}(\mathfrak{pr})$ can be different from $\mathfrak{s}$, \cf\ \cite[p. 2]{020R0}.

\begin{defi}[{\cite[Def. 4.5]{023BHKLS0}}]
The \emph{statistically refining\footnote{or statistically reaping} number} $\mathfrak{r}_*$ is the minimal cardinality of a statistically refining family.
\end{defi}

\begin{defi}[{\cite[p. 115]{984D0}, \cite[3.1 Def.]{010B0}, \cite[p. 182]{012H1}}]
\label{splitting number}
The \emph{splitting number} $\mathfrak{s}$ is the minimal cardinality of a splitting family.
\end{defi}

\begin{defi}[{\cite[p. 9]{001C0}}]
The $\QQ$-\emph{splitting number} $\mathfrak{s}_\QQ$ is the minimal cardinality of a $\QQ$-splitting family.
\end{defi}

\begin{defi}[{\cite[Introduction]{010HMM0}}]
The \emph{pair-splitting number} $\spair$ is the minimal cardinality of a pair-splitting family.
\end{defi}

\begin{defi}[{\cite[Definition 2.3]{023BHKLS0}}]
The \emph{bisecting number} $\mathfrak{s}_{\sfrac{1}{2}}$ is the minimal cardinality of a bisecting family.
\end{defi}

\begin{defi}[{\cite[Definition 4.4]{023BHKLS0}}]
The \emph{unbisecting (or semirefining\footnote{or semireaping}) number} $\mathfrak{r}_{\sfrac{1}{2}}$ is the minimal cardinality of an unbisected family.
\end{defi}

\begin{defi}[{\cite[Definition 2.3]{023BHKLS0}}]
The \emph{cofinally bisecting number} $\shalfinf$ is the minimal cardinality of a cofinally bisecting family.
\end{defi}

\begin{defi}[{\cite[Definition 2.3]{023BHKLS0}}]
The $\varepsilon$\emph{-almost bisecting number} $\shalfpme$ is the minimal cardinality of an $\varepsilon$-almost bisecting family.
\end{defi}

Recall that a series of real numbers $\displaystyle\sum_{i \nearrow \infty} r_i$ is \emph{conditionally convergent} if it is convergent but there is a set $A$ of natural numbers such that $\displaystyle\sum_{i \in A} r_i$ fails to be. Moreover we say that $\displaystyle\sum_{i \nearrow \infty} r_i$ \emph{diverges by oscillation} if for every natural number $k$ and every real number $r$ there are $\ell, m, n \in \omega \setminus k$ such that
\begin{align}
& \displaystyle\sum_{i < \ell} r_i \notin ]r - \sfrac{1}{k}, r + \sfrac{1}{k}[, \\
& \displaystyle\sum_{i < m} r_i > -k, \\
& \displaystyle\sum_{i < n} r_i < k.
\end{align}

\begin{defi}[{\cite[Definition 1]{019BBH0}}]
The \emph{subseries number} $\sss$ is the minimal cardinality of a family $\mathcal{I}$ of sets of natural numbers such that for every conditionally convergent series $\displaystyle\sum_{i \nearrow \infty} r_i$ there is an index set $I \in \mathcal{I}$ such that $\displaystyle\sum_{j \in I} r_j$ fails to converge to a real number.
\end{defi}

\begin{defi}[{\cite[Def. 2]{019BBH0}}]
The \emph{oscillating subseries number} $\sss_o$ is the minimal cardinality of a family $\mathcal{I}$ of sets of natural numbers such that for every conditionally convergent series $\displaystyle\sum_{i \nearrow \infty} r_i$ there is an index set $I \in \mathcal{I}$ such that $\displaystyle\sum_{j \in I} r_j$ diverges by oscillation.
\end{defi}

We call a pair $\pi = \opair{D}{\seq{\pi_n}{n \in D}}$ where $D$ is an infinite set of natural numbers and $\pi_n : \soafft{n}{\omega} \longrightarrow \omega$ for every $n \in D$, a \emph{predictor}. We say that $\pi$ \emph{predicts} a sequence $s \in \soafft{\omega}{\omega}$ if $\soast{n \in D}{\pi_n(s \upharpoonright n) \ne s(n)}$ is finite, otherwise $s$ \emph{evades} $\pi$.

\begin{defi}[{\cite[p. 529]{994B0}, \cite[10.1 Def.]{010B0}}]
\label{evasion number}
The \emph{evasion number} $\mathfrak{e}$ is the minimal cardinality of a family $\mathcal{E}$ of functions in $\soafft{\omega}{\omega}$ such that every predictor is evaded by a member of $\mathcal{E}$.
\end{defi}

\begin{defi}[{\cite[p. 529]{995B2}, \cite[10.2 Def.]{010B0}}]
\label{unbounded evasion number}
For a function $f : \omega \longrightarrow \omega \setminus 2$, the characteristic $\mathfrak{e}_f$ is the smallest cardinality of any family
\begin{align*}
\mathcal{E} \subset \prod_{n \in \omega} f(n)
\end{align*}
such that no single predictor predicts all members of $\mathcal{E}$. The \emph{unbounded evasion number} $\eubd$ is the minimum of $\mathfrak{e}_f$ over all functions $f : \omega \longrightarrow \omega \setminus 2$.
\end{defi}

\begin{defi}[{\cite[p. 57, Def. 16]{006MS0}, \cite[p. 97, Thm. 2.3]{981M0}}]
$\eqfr$ is the minimal cardinality of a bounded family $\mathcal{F} \subset \soafft{\omega}{\omega}$ such that for every $g \in \soafft{\omega}{\omega}$ there is an $f \in \mathcal{F}$ such that $f(n) \ne g(n)$ for all natural numbers $n$.
\end{defi}

\begin{defi}
Given a proper ideal $\mathcal{I}$ on a set $X$,
\begin{align*}
\add(\mathcal{I}) & \text{ is the minimal cardinality of a subfamily of } \mathcal{I} \text{ whose union is not in } \mathcal{I}, \\
\cov(\mathcal{I}) & \text{ is the minimal cardinality of a subfamily of } \mathcal{I} \text{ whose union is } X, \\
\non(\mathcal{I}) & \text{ is the minimal cardinality of a an element of } \Pot(X) \setminus \mathcal{I}, \\
\cof(\mathcal{I}) & \text{ is the minimal cardinality of a base of } \mathcal{I}.
\end{align*}

We have $\add(\mathcal{I}) \leqslant \cov(\mathcal{I})$. If $\mathcal{I}$ contains all singletons then, clearly, $\add(\mathcal{I}) \leqslant \non(\mathcal{I})$, also, the union of any base is then $X$, hence $\cov(\mathcal{I}) \leqslant \cof(\mathcal{I})$. Also, $\non(\mathcal{I}) \leqslant \cof(\mathcal{I})$ for otherwise one could inductively in step $\alpha$ pick an $x_\alpha \in X \setminus (B_\alpha \cup \soast{x_\beta}{\beta < \alpha})$ for any base $\soast{B_\alpha}{\alpha < \cof(\mathcal{I})} \in [\mathcal{I}]^{\cof(\mathcal{I})}$ which would yield $\soast{x_\alpha}{\alpha < \cof(\mathcal{I})} \in [X]^{\cof(\mathcal{I})} \setminus \mathcal{I}$. 
\end{defi}

\begin{defi}
\label{starred characteristics}
Given an ideal $\mathcal{I}$ on the set of natural numbers,
\begin{align*}
\add^*(\mathcal{I}) & \text{ is the minimal size of an } \mathcal{F} \subset \mathcal{I} \\
& \text{ such that for all } x \in \mathcal{I} \text{ there is a } y \in \mathcal{F} \text{ with } \card{x \setminus y} = \aleph_0.\\
\non^*(\mathcal{I}) & \text{ is the minimal size of an } \mathcal{F} \subset [\omega]^\omega \\
& \text{ such that for all } x \in \mathcal{I} \text{ there is a } y \in \mathcal{F} \text{ with } \card{x \cap y} < \aleph_0.\\
\cov^*(\mathcal{I}) & \text{ is the minimal size of an } \mathcal{F} \subset \mathcal{I} \\
& \text{ such that for all } x \in [\omega]^\omega \text{ there is a } y \in \mathcal{F} \text{ with } \card{y \setminus x} < \aleph_0.\\
\cof^*(\mathcal{I}) & \text{ is the minimal size of an } \mathcal{F} \subset \mathcal{I} \\
& \text{ such that for all } x \in \mathcal{I} \text{ there is a } y \in \mathcal{F} \text{ with } \card{x \cap y} < \aleph_0.
\end{align*}
\end{defi}

We now consider partitions of $\omega$. We say that a partition $R$ \emph{almost refines} another partition $X$ if
\begin{align*}
\soast{A \in X}{\forall B \subset R : \bigcup B \ne A} \text{ is finite.}
\end{align*}
We say that $S$ \emph{dual-splits} $X$ if $S$ is not almost refining $X$ yet there is an infinite partition of the natural numbers almost refining both $S$ and $X$.

\begin{defi}[{\cite[p. 28]{992CMW0}}]
The \emph{dual reaping number} $\mathfrak{r}_d$ is the minimal cardinality of a family $\mathcal{F}$ of partitions of $\omega$ for which there is no single partition of $\omega$ dual-splitting all members of $\mathcal{F}$.
\end{defi}

A family of sets is called \emph{linear}, if it is linearly ordered by $\subsetmf$ and a family $\mathcal{G}$ is a \emph{refinement} of a family $\mathcal{F}$ if for every $A \in \mathcal{F}$ there is a $B \in \mathcal{G}$ with $B \subsetmf A$.

\begin{defi}[{\cite[Def. 1.1]{016MST0}}]
The \emph{linear refinement number} $\mathfrak{lr}$ is the minimal cardinality of a family $\mathcal{F} \subset [\omega]^\omega$ without a linear refinement such that for every finite $\mathcal{A} \subset \mathcal{F}$ the set $\bigcap \mathcal{A}$ is infinite.
\end{defi}

\begin{defi}[{\cite[Def. 1.2]{016MST0}}]
The \emph{linear excluded middle number} $\mathfrak{lx}$ is the minimal cardinality of a family $\mathcal{F} \subset \soafft{\omega}{\omega}$ such that for every $h \in \soafft{\omega}{\omega}$ for which $\mathcal{S} \dfeq \bsoast{\soast{n < \omega}{f(n) \leqslant h(n)}}{f \in \mathcal{F}} \subset [\omega]^\omega$, the family $\mathcal{S}$ has no linear refinement.
\end{defi}

\subsubsection{A Summary and a Diagram}

\label{summary}
\begin{table}
$
\begin{array}{|rcl|rl|}
\hline
\eqfr & \leqslant & \displaystyle \inf_{\varepsilon \nearrow \sfrac{1}{2}}(\shalfpme) & \text{\upshape Theorem} & \ref{eqfr <= inf s_1/2+-e} \\
\hline
\mathfrak{r}_* & \leqslant & \cov(\smz) & \text{\upshape Theorem} & \ref{r_* <= cov(SN)} \\
\hline
\min(\mathfrak{b}, \mathfrak{r}_{\sfrac{1}{2}}) & \leqslant & \sss & \text{\upshape Theorem} & \ref{min(b, r_1/2) <= ß} \\
\hline
\mathfrak{s}(\mathfrak{pr}) & \leqslant & \non(\meagre) & \text{\upshape Theorem} & \ref{spr <= non[meagre]} \\
\hline
\mathfrak{e} & \leqslant & \spair & \text{\upshape Theorem} & \ref{e_ubd <= s_pair} \\
\hline
\spair & \leqslant & \sss & \text{\upshape Theorem} & \ref{s_pair <= ß} \\
\hline
\end{array}
$
\caption{Summary of Proven Results}
\end{table}

Note that the Mathias model witnesses both the consistency of $\mathfrak{h} = \aleph_2$, \cite[Section 11.8]{010B0}, and the one of the Borel conjecture, \cite[Thm. 8.3.4]{995BJ0}. So in particular, the Mathias model believes in $\eqfr = \non(\smz) = \aleph_1$ thus witnessing that the inequalities in Theorems \ref{eqfr <= inf s_1/2+-e}, and \ref{e_ubd <= s_pair} are consistently strict. Analogously to the proof of Theorem \ref{spr <= non[meagre]} and even more straightforwardly it is possible to show that $\mathfrak{s}(\mathfrak{pr}) \leqslant \mathfrak{d}$. See \cite[3.3 Thm.]{010B0} and \cite[Thm. 2.2]{020CGW0} to get an idea about proving it. Furthermore, the Random model witnesses that $\max(\mathfrak{b}, \mathfrak{d}, \non(\lmz)) < \cov(\lmz)$ is consistent, \cite[Section 11.4]{010B0}, so \emph{a fortiori} the inequalities in Theorems \ref{min(b, r_1/2) <= ß}, \ref{spr <= non[meagre]}, and \ref{s_pair <= ß} are consistently strict.

\begin{table}
\arraycolsep = 0dd
$\begin{array}{|rcccccl|c|c|}
\hline
&&&& \multicolumn{2}{c}{\headline{Inequalities}} && \headline{Origin} & \headline{Proof} \\
\hline
\cov(\lmz) \quad & \leqslant & \, \sss \, & \leqslant & \sss_o & \leqslant & \non(\meagre) & \text{\cite{019BBH0}} & \\
\hline
\cov(\lmz) \quad & \leqslant &\, \mathfrak{r}_{\sfrac{1}{2}} \, & \leqslant & \mathfrak{r}_* & \leqslant & \non(\meagre) & \text{\cite{023BHKLS0}} & \\
\hline
\mathfrak{s} \quad & \leqslant & \, \shalfpme \, & \leqslant & \mathfrak{s}_{\sfrac{1}{2}} & \leqslant & \non(\lmz) & \text{\cite{023BHKLS0}} & \\
\hline
\mathfrak{s} \quad & \leqslant & \, \swhalf\, & \leqslant & \shalfinf & \leqslant & \non(\lmz) & \text{\cite{023BHKLS0}} & \\
\hline
&&\mathfrak{s} & \leqslant & \spair & \leqslant & \non(\lmz) & \text{\cite{010M0}} & \\
\hline
&&&& \mathfrak{r}_d & \leqslant & \non(\meagre), \, \non(\lmz) & \text{J. Brendle} & \text{\cite[p. 503, Prop. 6]{010M0}} \\
\hline
&&&&\mathfrak{s}_\QQ & \leqslant & \add(\meagre), \, \mathfrak{s} & \text{\cite[p. 57]{006B2}} & \\
\hline
&&&&\mathfrak{h}_\QQ & \leqslant & \mathfrak{s}_\QQ & \text{\cite[p. 73]{004BHH0}} & \\
\hline
&&&&\mathfrak{e} & \leqslant & \cov(\meagre) & \text{\cite{998K0}} & \\
\hline
&&&&\mathfrak{h} & \leqslant & \mathfrak{g} & \text{\cite[p. 53]{989BL0}} & \text{\cite[p. 433, 6.27 Prop.]{010B0}} \\
\hline
&&&&\mathfrak{h} & \leqslant & \mathfrak{s} & \text{\cite[p. 354]{989BS0}} & \text{\cite[p. 427, 6.9 Thm.]{010B0}} \\
\hline
&&&\multicolumn{2}{c}{\cov(\meagre)} & \leqslant & \eqfr & \text{\cite{987B1}} & \text{\cite[p. 420, 5.9 Thm.]{010B0}} \\
\hline
&&&&\mathfrak{b} & \leqslant & \non(\meagre) & \text{\cite[p. 237]{987B1}} & \text{\cite[p. 419, 5.5 Prop.]{010B0}} \\
\hline
&&&&\mathfrak{s} & \leqslant & \non(\meagre) & \text{Folklore} & \text{\cite[p. 424, 5.19 Thm.]{010B0}} \\
\hline
&&\multicolumn{3}{c}{\max(\mathfrak{b}, \mathfrak{s}, \mathfrak{lr})} & \leqslant & \mathfrak{lx} & \text{\cite[Cor. 2.13]{016MST0}} & \\
\hline
&&&&\cov^*(\mathcal{Z}_0) & \leqslant & \max\big(\mathfrak{b}, \mathfrak{s}(\mathfrak{pr})\big) & \text{\cite[Cor. 30]{020R0}} & \\
\hline
\end{array}$
\label{legend}
\caption{Legend}
\end{table}

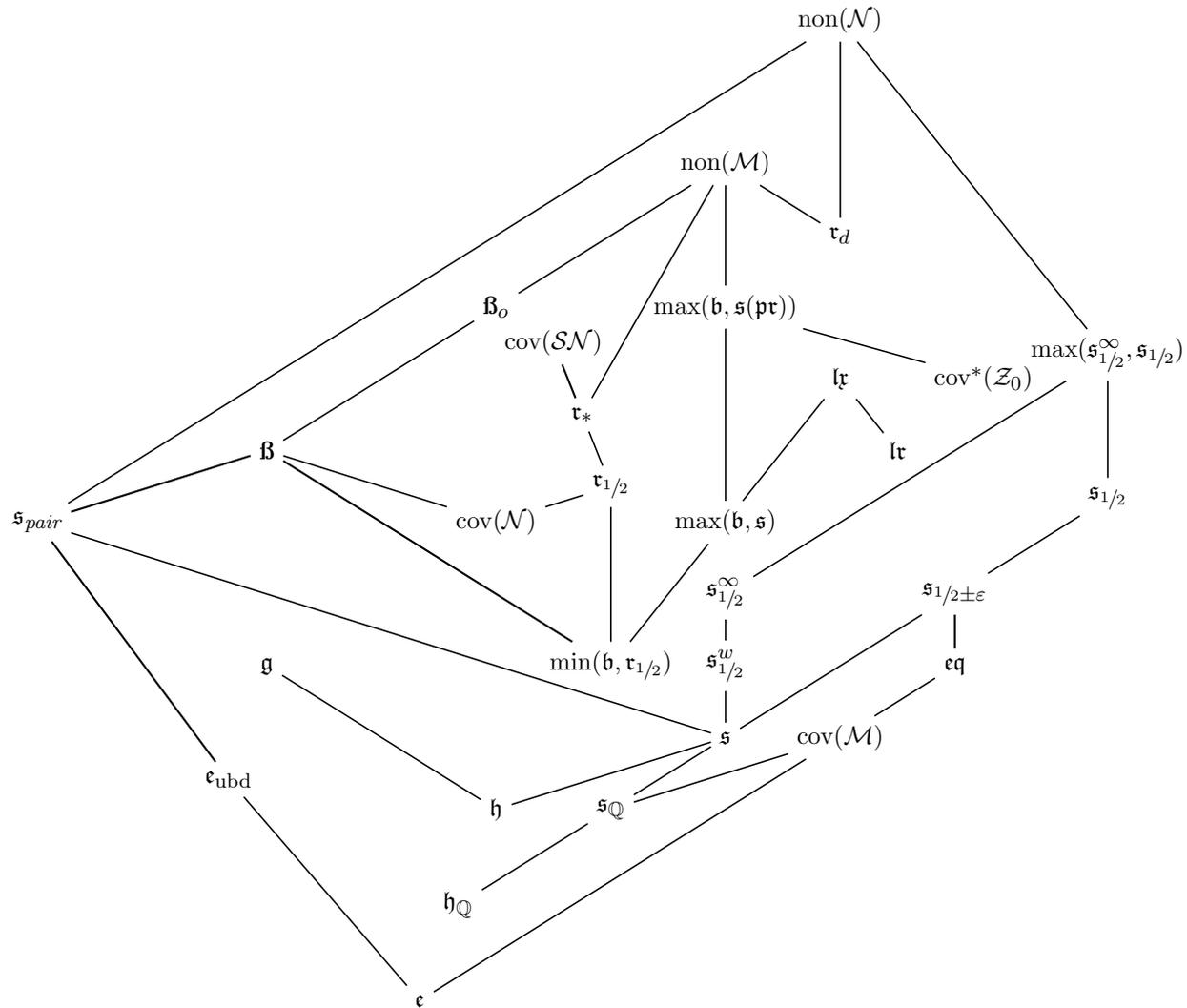
\begin{figure}
\label{diagram}
\begin{center}
\footnotesize
\begin{tikzpicture}[line/.style={semithick}, proven/.style={thick}]
\pgftransformxscale{0.8}
\pgftransformyscale{0.5}
\node (ss) at (-10, 6) {$\sss$};
\node (sso) at (-6, 10) {$\sss_o$};
\node (s) at (-2, -2) {$\mathfrak{s}$};
\node (sq) at (-4, -4) {$\mathfrak{s}_\QQ$};
\node (s12inf) at (-2, 2) {$\shalfinf$};
\node (non[L]) at (0, 18) {$\non(\lmz)$};
\node (mcblb) at (4.667, 8.667) {$\max(\shalfinf, \mathfrak{s}_{\sfrac{1}{2}})$};
\node (fs) at (-2, 4) {$\max(\mathfrak{b}, \mathfrak{s})$};
\node (max[bspr]) at (-2, 10) {$\max(\mathfrak{b}, \mathfrak{s}(\mathfrak{pr}))$};
\node (cov[Z]) at (2.5, 8) {$\cov^*(\mathcal{Z}_0)$};
\node (lx) at (0, 8) {$\mathfrak{lx}$};
\node (lr) at (1, 6) {$\mathfrak{lr}$};
\node (non[M]) at (-2, 14) {$\non(\meagre)$};
\node (cov[SN]) at (-5, 9) {$\cov(\smz)$};
\node (r_d) at (0, 12) {$\mathfrak{r}_d$};
\node (cov[L]) at (-6, 4) {$\cov(\lmz)$};
\node (cov[M]) at (0, -2) {$\cov(\meagre)$};
\node (h) at (-6, -4) {$\mathfrak{h}$};
\node (g) at (-10, 0) {$\mathfrak{g}$};
\node (hq) at (-6.667, -6.667) {$\mathfrak{h}_\QQ$};
\node (eq) at (2, 0) {$\eqfr$};
\node (s12+-e) at (2, 2) {$\shalfpme$};
\node (s12) at (4.667, 4.667) {$\mathfrak{s}_{\sfrac{1}{2}}$};
\node (s12w) at (-2, 0) {$\swhalf$};
\node (min[br12]) at (-4, 0) {$\min(\mathfrak{b}, \mathfrak{r}_{\sfrac{1}{2}})$};
\node (r_1/2) at (-4, 5) {$\mathfrak{r}_{\sfrac{1}{2}}$};
\node (r*) at (-4.5, 7) {$\mathfrak{r}_*$};
\node (e) at (-7.333, -9.333) {$\mathfrak{e}$};
\node (eu) at (-10.667, -3.167) {$\eubd$};
\node (sp) at (-14, 4) {$\spair$};
\draw (sp) [line] -- (non[L]);
\draw (r_d) [line] -- (non[M]);
\draw (r_d) [line] -- (non[L]);
\draw (mcblb) [line] -- (non[L]);
\draw (s12) [line] -- (mcblb);
\draw (s12inf) [line] -- (mcblb);
\draw (cov[L]) [line] -- (ss);
\draw (cov[L]) [line] -- (r_1/2);
\draw (min[br12]) [line] -- (r_1/2);
\draw (min[br12]) [line] -- (fs);
\draw (s) [line] -- (s12+-e);
\draw (s12+-e) [line] -- (s12);
\draw (s) [line] -- (sp);
\draw (s) [line] -- (s12w);
\draw (s12w) [line] -- (s12inf);
\draw (sq) [line] -- (s);
\draw (ss) [line] -- (sso);
\draw (sso) [line] -- (non[M]);
\draw (r_1/2) [line] -- (r*);
\draw (r*) [line] -- (non[M]);
\draw (max[bspr]) [line] -- (non[M]);
\draw (fs) [line] -- (max[bspr]);
\draw (fs) [line] -- (lx);
\draw (lr) [line] -- (lx);
\draw (cov[Z]) [line] -- (max[bspr]);
\draw (cov[M]) [line] -- (eq);
\draw (e) [line] -- (cov[M]);
\draw (sq) [line] -- (cov[M]);
\draw (hq) [line] -- (sq);
\draw (h) [line] -- (g);
\draw (h) [line] -- (s);
\draw (min[br12]) [proven] -- (ss);
\draw (r*) [proven] -- (cov[SN]);
\draw (sp) [proven] -- (ss);
\draw (eq) [proven] -- (s12+-e);
\draw (e) [line] -- (eu);
\draw (eu) [proven] -- (sp);
\end{tikzpicture}
\end{center}
\caption{A diagram of cardinal characteristics}
\end{figure}

\clearpage
\section{Results}

\subsection{Inequalities Involving Halfway New Cardinal Characteristics}

\subsubsection{Splitting}

\begin{theo}
\label{eqfr <= inf s_1/2+-e}
$\displaystyle\eqfr \leqslant \inf_{\varepsilon \nearrow \sfrac{1}{2}}(\shalfpme)$.
\end{theo}

The following proof uses a partition of the natural numbers into intervals whose lengths are growing sufficiently quickly so that one can define an eventually different sequence from a bisecting sequence.

\begin{proof}
Suppose towards a contradiction that $\shalfpme < \eqfr$ for some $\varepsilon \in {]0, \sfrac{1}{2}[}$ and let $\mathcal{S}$ be an $\varepsilon$-almost bisecting family of cardinality $\shalfpme$. Let $\seq{I_k}{k < \omega}$ be the interval partition of the natural numbers defined by $\min(I_0) = 0$, $\min(I_{k + 1}) = \max(I_k) + 1$ and $\max(I_k) = k + 2\lceil 2\min(I_k) / (1 - 2 \varepsilon)\rceil$ for all natural numbers $k$. Note that every interval has an odd length. For every natural number $k$, choose a function $h_k : \Pot(I_k) \longrightarrow 2^{\max(I_k) - \min(I_k)}$ such that $h_k(A) = h_k(B)$ \ifff\ $A \in \{B, I_k \setminus B\}$. Now we can define
\begin{align*}
\mathcal{E} \dfeq \soast{e \in \soafft{\omega}{\omega}}{\exists S \in \mathcal{S} \, \exists k \in \omega \, \forall \ell \in \omega \setminus k : e(\ell) = h_\ell(S \cap I_\ell)}.
\end{align*}
We have $\card{\mathcal{E}} = \aleph_0 \cdot \card{\mathcal{S}} = \aleph_0 \cdot \shalfpme = \shalfpme < \eqfr$. As $\mathcal{E}$ is bounded, by definition of $\eqfr$, there is an $f \in \soafft{\omega}{\omega}$ such that for all $e \in \mathcal{E}$ there is a natural number $\ell$ such that $e(\ell) = f(\ell)$. Now, for every natural number $k$, let $g_k$ be the right-inverse 
of $h_k$ such that $\card{g_k(\ell)} > \card{I_k \setminus g_k(\ell)}$ for all $\ell < 2^{\max(I_k) - \min(I_k)}$. We consider the set
\begin{align*}
X & \dfeq \bigcup_{\ell < \omega} g_\ell(\max(f(\ell), 2^{\max(I_\ell) - \min(I_\ell)} - 1)).
\intertext{Now there is an $S \in \mathcal{S}$ which is $\varepsilon$-almost bisecting $X$. Let $m$ be a natural number such that}
\frac{\card{S \cap X \cap n}}{\card{X \cap n}} & \in \left]\frac{1}{2} - \varepsilon, \frac{1}{2} + \varepsilon\right[ \text{ for all } n \in \omega \setminus m.
\intertext{Let $k$ be the least natural number such that $I_k \subset \omega \setminus m$. We define a function,}
d : \omega & \longrightarrow \omega \\
n & \longmapsto \begin{cases}
f(n) + 1 & \text{if } n < k, \\
h_n(S \cap I_n) & \text{else.}
\end{cases}
\intertext{Clearly, $d \in \mathcal{E}$. Now there is a natural number $\ell$ such that $d(\ell) = f(\ell)$. By definition of $d$ we have $\ell \in \omega \setminus k$, moreover}
f(\ell) & =  {d(\ell) = h_\ell(S \cap I_\ell) < 2^{\max(I_\ell) - \min(I_\ell)}}.
\intertext{Next we distinguish two cases. First assume that $2\card{S \cap I_\ell} < \card{I_\ell}$. Then}
X \cap I_\ell & = g_\ell(\max(f(\ell), 2^{\max(I_\ell) - \min(I_\ell)} - 1)) = g_\ell(f(\ell)) \\
& = g_\ell(d(\ell)) = g_\ell(h_\ell(I_\ell \cap S)) = I_\ell \setminus S, \\
& \text{and therefore } S \cap X \cap \min(I_{\ell + 1}) \subset \min(I_\ell) \text{. We have}
\end{align*}
\begin{align*}
\frac{\card{X \cap S \cap \min(I_{\ell + 1})}}{\card{X \cap \min(I_{\ell + 1})}} = \frac{\card{X \cap S \cap \min(I_\ell)}}{\card{X \cap \min(I_{\ell + 1})}} & \leqslant \frac{\min(I_\ell)}{\card{X \cap \min(I_{\ell + 1})}} \\
& < \frac{2\min(I_\ell)}{\min(I_{\ell + 1})} <  \frac{2\min(I_\ell)}{\max(I_\ell)} \leqslant \frac{2\min(I_\ell)(1 - 2\varepsilon)}{\ell(1 - 2\varepsilon) + 4\min(I_\ell)} \leqslant \frac{1}{2} - \varepsilon, \\
& \text{implying } \ell < k \text{, a contradiction.}
\end{align*}
Now we assume $\card{I_\ell} < 2\card{S \cap I_\ell}$, we have
\begin{align*}
X \cap I_\ell & = g_\ell(\max(f(\ell), 2^{\max(I_\ell) - \min(I_\ell)} - 1)) = g_\ell(f(\ell)) \\
= & g_\ell(d(\ell)) = g_\ell(h_\ell(I_\ell \cap S)) = I_\ell \cap S, \\
& \text{and therefore } S \cap I_\ell \subset X \cap S \text{. Then} \\
\frac{\card{X \cap S \cap \min(I_{\ell + 1})}}{\card{X \cap \min(I_{\ell + 1})}} & \geqslant \frac{\card{S \cap I_\ell}}{\card{X \cap \min(I_{\ell + 1})}} \geqslant \frac{\card{S \cap I_\ell}}{\min(I_\ell) + \card{X \cap I_\ell}} \\
& = \frac{\card{S \cap I_\ell}}{\min(I_\ell) + \card{S \cap I_\ell}} \text{ and as for nonnegative values,} \\
& \text{the function } \frac{x}{a + x} \text{ is monotonically increasing in } x, \\
& \geqslant \frac{\card{I_\ell}}{2\min(I_\ell) + \card{I_\ell}} = \frac{\min(I_{\ell + 1}) - \min(I_\ell)}{\min(I_{\ell + 1}) + \min(I_\ell)} \\
& > \frac{\max(I_\ell) - \min(I_\ell)}{\max(I_\ell) + \min(I_\ell)} \geqslant \frac{4\min(I_\ell) + (1 - 2\varepsilon)(\ell - \min(I_\ell))}{4\min(I_\ell) + (1 - 2\varepsilon)(\ell + \min(I_\ell))} \\
& \text{and as, for nonnegative values, the function } \frac{a + b(x - 1)}{a + b(x + 1)} \\
& \text{is monotonically increasing in } x,\\
& \geqslant \frac{3 + 2\varepsilon}{5 - 2\varepsilon} > \frac{1}{2} + \varepsilon \text{, implying } \ell < k \text{, a contradiction.}
\end{align*}
\end{proof}

\subsubsection{Reaping}

The following Lemma is strongly inspired by \cite[Lemma 8.1.13]{995BJ0}.

\begin{lemm}
\label{smz}
If $Z$ is a strong measure zero set in the Cantor space, $\seq{J_i}{i < \omega}$ is an interval partition of $\omega$, and $\soast{P_i}{i < \omega}$ is a family of infinitely many disjoint infinite sets of natural numbers, then there is a sequence $\seq{s_j}{j < \omega}$ such that for all natural numbers $k$ and all $x \in Z$ there is a $j \in P_k$ such that
\begin{align}
x \upharpoonright \bigcup_{i \leqslant j} J_i = s_j.
\end{align}
\end{lemm}

\begin{proof}
We consider the Cantor-space with the usual metric
\begin{align}
d : \soafft{\omega}{2} \times \soafft{\omega}{2} & \longrightarrow \RR, \\
\opair{x}{y} & \longmapsto \min(\soast{n < \omega}{x(n) \ne y(n)}).
\end{align}
By definition of the metric $d$, an interval of length $2^{-n}$ consists of all elements of the Cantor space sharing a certain sequence of $n$ bits as an initial segment. Now, for every natural number $n$, let $\seq{k_\opair{n}{i}}{i < \omega}$ be an enumeration of
\begin{align}
\soast{\bigcup_{i \leqslant m} J_i}{m \in P_n}.
\end{align}
Note that as $\seq{J_i}{i < \omega}$ is an interval partition, these $k_\opair{n}{i}$ are natural numbers. Moreover, as the $P_i$'s are disjoint, $k_\opair{n}{i} = k_\opair{m}{j}$ implies $m = n$ and $i = j$. Now consider the sequences $\seq{\varepsilon_\opair{n}{i}}{i < \omega}$ given by $\varepsilon_\opair{n}{i} = 2^{-k_\opair{n}{i}}$. By the definition of strong measure zero sets there are sequences $\seq{r_\opair{n}{i}}{i < \omega}$ of initial segments such that for all natural numbers $n$ and $i$, the segment $r_\opair{n}{i}$ is $k_\opair{n}{i}$ sequents long and for every $x \in Z$ and every natural number $n$ there is an $i < \omega$ such that $r_\opair{n}{i}$ is an initial segment of $x$. Now simply let $\seq{\ell_i}{i < \omega}$ be the increasing enumeration of $\soast{k_\opair{n}{i}}{\opair{n}{i} \in \omega \times \omega}$ and for every natural number $j$, let $s_j = r_\opair{n_j}{i_j}$ where $n_j$ and $i_j$ are such that $\ell_j = k_\opair{n_j}{i_j}$. The sequence $\seq{s_j}{j < \omega}$ provides what was demanded.
\end{proof}

\begin{theo}
\label{r_* <= cov(SN)}
$\mathfrak{r}_* \leqslant \cov(\smz)$.
\end{theo}

The following proof uses the fact that a bisecting real can be used to code avoidance of a sequence of sufficiently small open intervals.

\begin{proof}
Assume towards a contradiction that $\cov(\smz) < \mathfrak{r}_*$ and let $\mathcal{C}$ be a family of $\cov(\smz)$ sets of strong measure zero covering the Cantor space. Let $\seq{J_i}{i < \omega}$ be an interval partition defined by $\displaystyle J_i \dfeq 1 + i + (i + 2)! \setminus \bigcup_{k < i} J_k$. Note that all parts of this partition are of odd length! Furthermore, let $\seq{P_i}{i < \omega}$ be a partition of $\omega$ into infinitely many infinite parts with $P_i \subset \omega \setminus i$ for all natural numbers $i$. For every $Z \in \mathcal{C}$, use \autoref{smz} to find a sequence $\seq{s^Z_j}{j < \omega}$ of $s^Z_j \in \soafft{\bigcup_{i \leqslant j} J_i}{2}$ such that for all natural numbers $k$ and all $x \in Z$ there is a $j \in P_k$ such that $x \upharpoonright \bigcup_{i \leqslant j} J_i = s^Z_j$. Now we can define
\begin{align}
\label{D-def} \mathcal{D} \dfeq & \soast{X_Z}{Z \in \mathcal{C}} \\
\text{where } X_Z \dfeq & \bigcup_{i < \omega} F^Z_i \\
\text{and } F^Z_i & \text{ is the larger of the two sets} \\
\soast{k \in J_i}{s^Z_i(k) = 0} & \text{ and } \soast{k \in J_i}{s^Z_i(k) = 1}.
\end{align}
As $\card{\mathcal{D}} \leqslant \card{\mathcal{C}} = \cov(\smz) < \mathfrak{r}_*$, there is a moderate $S \in [\omega]^\omega$ such that
\begin{align}
\label{lim}
\displaystyle 1 = & \lim_{n \nearrow \infty} \frac{n \cdot \card{S \cap X \cap n}}{\card{S \cap n} \cdot \card{X \cap n}} \text{ for all } X \in \mathcal{D}.
\intertext{
As $\mathcal{C}$ covers the Cantor space, we can fix a $Z \in \mathcal{C}$ containing $\chi_S$. Let $\varepsilon \in {]0, \sfrac{1}{2}[}$ be such that
}
\varepsilon < & \liminf_{n \nearrow \infty} \frac{\card{S \cap n}}{n} \text{ and } \limsup_{n \nearrow \infty} \frac{\card{S \cap n}}{n} < 1 - \varepsilon.
\intertext{As $Z \in \mathcal{C}$, by \eqref{D-def} we have $X_Z \in \mathcal{D}$ and by 
\eqref{lim} we then also have }
\varepsilon < & \liminf_{n \nearrow \infty} \frac{\card{S \cap X_Z \cap n}}{\card{X_Z \cap n}} \text{ and }\displaystyle\limsup_{n \nearrow \infty} \frac{\card{S \cap X_Z \cap n}}{\card{X_Z \cap n}} < 1 - \varepsilon.
\end{align}
Now choose a natural number $k > \sfrac{3}{\varepsilon}$ large enough such that $\displaystyle\frac{\card{S \cap X_Z \cap n}}{\card{X_Z \cap n}} \in \left]\varepsilon, 1 - \varepsilon\right[$ for all $n \in \omega \setminus k$. Choose a natural number $j$ such that $P_j \subset \omega \setminus k$ and let $\ell \in P_j$ be such that $\chi_S \upharpoonright \bigcup_{i \leqslant \ell} J_i = s^Z_\ell$. Note that $\displaystyle 2\card{X_Z \cap \bigcup_{i \leqslant \ell} J_i} > \card{\bigcup_{i \leqslant \ell} J_i}$. Also note that, as $P_j \ni \ell$ consists of natural numbers no smaller than $k$, we have $\ell \geqslant k > \sfrac{3}{\varepsilon}$. We now distinguish two cases. First assume, that $F^Z_\ell = \soast{k \in J_\ell}{s^Z_\ell(k) = 0}$. Then $S \cap X_Z \cap \bigcup_{i \leqslant \ell} J_i \subset \bigcup_{i < \ell} J_i$ and hence
\begin{align*}
\frac{\card{S \cap X_Z \cap \bigcup_{i \leqslant \ell} J_i}}{\card{X_Z \cap \bigcup_{i \leqslant \ell} J_i}} & \leqslant \frac{2\card{\bigcup_{i < \ell} J_i}}{\card{\bigcup_{i \leqslant \ell} J_i}}
= \frac{2(\ell + (\ell + 1)!)}{1 + \ell + (\ell + 2)!} < \frac{2(\ell + (\ell + 1)!)}{\ell + (\ell + 2)!} \\
& = \frac{2(1+ (\ell + 1)(\ell - 1)!)}{1 + (\ell + 1)(\ell + 2)(\ell - 1)!} < \frac{2(1+ (\ell + 1)(\ell - 1)!)}{(\ell + 1)(\ell + 2)(\ell - 1)!} \\
& \leqslant \frac{2(\ell + 2)(\ell - 1)!}{(\ell + 1)(\ell + 2)(\ell - 1)!} = \frac{2}{\ell + 1} < \frac{2}{\ell} < \frac{2}{\sfrac{3}{\varepsilon}} = \frac{2\varepsilon}{3} < \varepsilon.
\intertext{Now assume, that $F^Z_\ell = \soast{k \in J_\ell}{s^Z_\ell(k) = 1}$. Then $X_Z \cap J_\ell \subset S \cap X_Z \cap \bigcup_{i \leqslant \ell} J_i$ and hence}
\frac{\card{S \cap X_Z \cap \bigcup_{i \leqslant \ell} J_i}}{\card{X_Z \cap \bigcup_{i \leqslant \ell} J_i}} &\geqslant \frac{\card{X_Z \cap J_\ell}}{\card{X_Z \cap \bigcup_{i \leqslant \ell} J_i}} = \frac{\card{X_Z \cap J_\ell}}{\card{X_Z \cap \bigcup_{i <\ell} J_i} + \card{X_Z \cap J_\ell}}. \\
\text{As } X_Z \cap J_\ell = F^Z_\ell & \text{ and } 2\card{F^Z_\ell} > \card{J_\ell} \text{ and for nonnegative values, the function} \\
& \frac{x}{a + x} \text{ is monotonically increasing in } x, \\
& \frac{\card{X_Z \cap J_\ell}}{\card{X_Z \cap \bigcup_{i <\ell} J_i} + \card{X_Z \cap J_\ell}} > \frac{\card{J_\ell}}{2\card{X_Z \cap \bigcup_{i <\ell} J_i} + \card{J_\ell}} \\
\geqslant & \frac{\card{J_\ell}}{2\card{\bigcup_{i <\ell} J_i} + \card{J_\ell}} = \frac{1 + (\ell + 2)! - (\ell + 1)!}{(\ell + 2)! + (\ell + 1)! + 2\ell + 1} \\
> & \frac{(\ell + 2)! - (\ell + 1)!}{(\ell + 2)! + 2(\ell + 1)!} = \frac{\ell + 1}{\ell + 4} > \frac{3 + \varepsilon}{3 + 4\varepsilon} > 1 - \varepsilon.
\end{align*}
In both cases we reach a contradiction thus concluding the proof.
\end{proof}

\begin{lemm}
\label{lemm : series}
\begin{align}
\text{The series } \sum_{n < \omega} a_n & \\
\text{with } a_n & = \frac{(-1)^n 2n}{(n + 1)(2n + 1)} \text{ is conditionally convergent.}
\end{align}
\end{lemm}

\begin{proof}
The sequence is convergent by the Leibniz criterion as it is alternating and its terms (ignoring $a_0$, which is $0$) decrease in absolute value. A failure of the latter would mean that
\begin{alignat}{2}
& \frac{2n}{(n + 1)(2n + 1)} & \leqslant & \quad \frac{2(n + 1)}{(n + 2)(2(n + 1) + 1))} \\
\Longrightarrow \quad & 2n(n + 2)(2n + 3) & \leqslant & \quad 2(n + 1)^2(2n + 1) \\
\Longrightarrow \quad & 4n^3 + 14n^2 + 12n & \leqslant & \quad 4n^3 + 10n^2 + 8n + 2 \\
\Longrightarrow \quad & 2n^2 + 2n & \leqslant & \quad 1 \\
\Longrightarrow \quad & n = 0. \\
\intertext{The series is, however, not absolutely convergent as} \\
 \absv{a_n} & \geqslant \frac{1}{2(n + 1)} \text{ for } n \in \omega \setminus 1. 
\end{alignat}
\end{proof}

\begin{lemm}
\label{lemm : permutation}
\begin{align}
\text{If } & \sum_{n < \omega} a_n \text{ is a series,} \\
\seq{I_n}{n < \omega} & \text{ is an interval partition of } \omega, \\
\text{and } b_n = & \sum_{i \in I_n} a_k \text{ is such that for all } n, \\
& \sum_{n < \omega} b_n \text{ is conditionally convergent,} \\
\text{then } & \sum_{n < \omega} a_n \text{ is conditionally convergent.}
\end{align}
\end{lemm}

\begin{proof}
Let $\sum_{n < \omega} a_n$, $\seq{I_n}{n < \omega}$, and $b_n (n < \omega)$ be as above. We may assume \wlg\ that $\sum_{n < \omega} b_n$ converges to a real number $r$. Let $\pi : \omega \longleftrightarrow \omega$ be a permutation such that
\begin{align}
\sum_{n < \omega} & b_{\pi(n)} = s \text{ where } s \ne r.
\intertext{Let $\rho : \omega \longleftrightarrow \omega$ be a permutation such that $\rho(n) \in I_{\pi(k)} \longleftrightarrow n \in I_k$. Then, clearly,}
\sum_{n < \omega} & a_{\rho(n)} = s \text{, while } \sum_{n < \omega} a_n = r \text{, so} \\
\sum_{n < \omega} & a_n \text{ is conditionally convergent.}
\end{align}
\end{proof}

The following lemma is a formulation of the technique of ``padding with zeroes" as introduced in \cite[Section 3]{020BBBHHL0}

\begin{lemm}
\label{padding}
\begin{align}
\text{If } \seq{n_k}{k < \omega} \text{ is an ascending sequence of natural numbers,} \\
\text{and } \sum_{n < \omega} a_n \text{ is conditionally convergent,} \\
\text{then } \sum_{n < \omega} b_n \text{ is conditionally convergent,} \\
\text{where } b_m = \begin{cases}
0 \text{ if } \forall k < \omega : m \ne n_k, \\
a_k \text{ if } m = n_k.
\end{cases}
\end{align}
\end{lemm}

\begin{proof}
Left to the reader.
\end{proof}


\begin{theo}
\label{min(b, r_1/2) <= ß}
$\min(\mathfrak{b}, \mathfrak{r}_{\sfrac{1}{2}}) \leqslant \sss$.
\end{theo}

The idea of the following proof is to use a bound for the distances between elements of the index set of a subseries and observe that a bisecting sequence implies convergence of the respective subseries. An earlier version of this preprint contained a version of this proof in which the $n^I_j$'s were defined using linearly shrinking intervals instead of quadratically shrinking ones. Due to a miscalculation this seemed to work.

\begin{proof}
Assume towards a contradiction that $\sss < \min(\mathfrak{b}, \mathfrak{r}_{\sfrac{1}{2}})$. Let $\mathcal{I}$ be a $\sss$-sized family of infinite index sets of natural numbers such that for every conditionally convergent series $\displaystyle\sum_{i \rightarrow \infty} r_i$ there is an index set $I \in \mathcal{I}$ such that $\displaystyle\sum_{i \in I} r_i$ is diverging. As $\card{\mathcal{I}} = \sss < \mathfrak{r}_{\sfrac{1}{2}}$, there is an $S \in [\omega]^\omega$ bisecting all members of $\mathcal{I}$. We define a family
\begin{align*}
\mathcal{B} & \dfeq \soast{f_I}{I \in \mathcal{I}}, \text{ where} \\
f_I & \dfeq \seq{n^I_j}{j < \omega} \text{, and} \\
n^I_j & \dfeq \min\left(\bgsoast{k < \omega}{\forall \ell \in \omega \setminus k : \frac{\card{S \cap I \cap \ell}}{\card{I \cap \ell}} \in \left]\frac{1}{2} - \frac{1}{4j + 12}, \frac{1}{2} + \frac{1}{4j + 12}\right[}\right).
\end{align*}
As $\card{\mathcal{B}} \leqslant \card{\mathcal{I}} = \sss < \mathfrak{b}$, there is a $g \in \soafft{\omega}{\omega}$ such that $f \leqmf g$ for all $f \in \mathcal{B}$. We may suppose \wlg\ that $g(i)2(i + 1) \leqslant g(i + 1)$ for all natural numbers $i$. We define a series
\begin{align*}
\sum_{k \rightarrow \infty} t_k & \\
\text{with } t_k & = \begin{cases}
\frac{(-1)^i}{g(i)(i + 1)(2i + 1)} \text{ for } k \in g(i)2(i + 1) \setminus g(i), \\
0 \text{ else.}
\end{cases}
\end{align*}
By Lemmata \ref{lemm : series}, \ref{lemm : permutation}, and \autoref{padding}, this series is conditionally convergent. By assumption, there is an index set $I \in \mathcal{I}$ such that $\displaystyle\sum_{j \in I} t_j$ is diverging. Let $k \in \omega \setminus 1$ be such that $f_I(\ell) \leqslant g(\ell)$ for all $\ell \in \omega \setminus k$.
Since
\begin{align*}
\sum_{j \in I} t_j = \sum_{j \nearrow \infty} u_j & \text{ with} \\
u_j & \dfeq \sum_{\ell \in I \cap g(j)2(j + 1) \setminus g(j)} t_\ell,
\end{align*}
and $u_j < 0$ \ifff\ $j$ is odd, the Leibniz criterion implies that there has to be a $j \in \omega \setminus k$ such that $\absv{u_j} \leqslant \absv{u_{j + 1}}$. But then
\begin{align*}
& \sum_{\ell \in I \cap g(j)2(j + 1) \setminus g(j)} \absv{t_\ell} = \absv{u_j} \leqslant \absv{u_{j + 1}} = \sum_{\ell \in I \cap g(j + 1)2(j + 2) \setminus g(j + 1)} \absv{t_\ell} \\
\Longleftrightarrow & \sum_{\ell \in I \cap g(j)2(j + 1) \setminus g(j)} \frac{1}{g(j)(j + 1)(2j + 1)} = \absv{u_j} \leqslant \absv{u_{j + 1}} \\
= & \sum_{\ell \in I \cap g(j + 1)2(j + 2) \setminus g(j + 1)} \frac{1}{g(j + 1)(j + 2)(2j + 3)} \\
\Longleftrightarrow & \frac{\card{I \cap g(j)2(j + 1) \setminus g(j)}}{g(j)(j + 1)(2j + 1)} = \absv{u_j} \leqslant \absv{u_{j + 1}} = \frac{\card{I \cap g(j + 1)2(j + 2) \setminus g(j + 1)}}{g(j + 1)(j + 2)(2j + 3)}.
\intertext{We have $k \leqslant j$ and therefore $f_I(j) \leqslant g(j)$, so by definition of $f_I$, the absolute values of $u_j$ and $u_{j + 1}$ can be estimated as follows:}
\absv{u_j} \geqslant & \, \frac{\left(\sfrac{1}{2} - \frac{1}{4j + 12}\right)g(j)2(j + 1) - \left(\sfrac{1}{2} + \frac{1}{4j + 12}\right)g(j)}{g(j)(j + 1)(2j + 1)} \\
 = & \frac{(2j + 5)2(j + 1) - (2j + 7)}{(4j + 12)(j + 1)(2j + 1)} = \frac{4j^2 + 12j + 3}{(4j + 12)(2j^2 + 3j + 1)} \\ 
\absv{u_{j + 1}} \leqslant & \, \frac{\left(\sfrac{1}{2} + \frac{1}{4j + 12}\right)g(j + 1)2(j + 2) - \left(\sfrac{1}{2} - \frac{1}{4j + 12}\right)g(j + 1)}{g(j + 1)(j + 2)(2j + 3)} \\
= & \frac{(2j + 7)2(j + 2) - (2j + 5)}{(4j + 12)(j + 2)(2j + 3)} = \frac{4j^2 +20j + 23}{(4j + 12)(2j^2 + 7j + 6)}.
\end{align*}
Therefore,
\begin{alignat*}{3}
& (4j + 12)(2j^2 + 3j + 1)(2j^2 + 7j + 6)\absv{u_j} && \leqslant (4j + 12)(2j^2 + 3j + 1)(2j^2 + 7j + 6)\absv{u_{j + 1}} \\
\Longrightarrow & (4j^2 + 12j + 3)(2j^2 + 7j + 6) && \leqslant (4j^2 + 20j + 23)(2j^2 + 3j + 1) \\
\Longleftrightarrow & 8j^4 + 52j^3 + 114j^2 + 93j + 18 && \leqslant 8j^4 + 52j^3 + 110j^2 + 89j + 23 \\
\Longleftrightarrow & 4j^2 + 4j && \leqslant 5 \\
\Longrightarrow & 4j(j + 1) && \leqslant 5 \text{, a contradiction.}
\end{alignat*}
\end{proof}

\subsection{Various Inequalities}

\subsubsection{Splitting of Partitions and the Uniformity of the Meagre Ideal}

\begin{theo}
\label{spr <= non[meagre]}
$\mathfrak{s}(\mathfrak{pr}) \leqslant \non(\meagre)$.
\end{theo}

The following proof uses Bartoszy{\'n}ski's characterisation of $\non(\meagre)$ via sequences of natural numbers equalling each other infinitely often. Such sequences are used in this proof to code sets of natural numbers thus guessing larger and larger parts of an infinite set of natural numbers thus enabling the definition of a splitting partition.

\begin{proof}
Assume towards a contradiction that $\non(\meagre) < \mathfrak{s}(\mathfrak{pr})$. Recall that $\non(\meagre)$ is the minimal cardinality of an infinitely often equal family, \cf\ \cite[5.9 Thm.]{010B0}. So let $\mathcal{E}$ be an infinitely often equal family of size $\non(\meagre)$. Let a bijection $e : \omega \longrightarrow [\omega]^{< \omega}$ encode all finite sets of natural numbers by a natural number. For every sequence $s \in \soafft{\omega}{\omega}$ we inductively define an infinite partition $P_s$ of $\omega$. If $i < j$ or $\card{e(s(i))} \ne (2i)!(i + 1)$ let $Q_s(j, i) \dfeq 0$, otherwise let $\seq{f_s(i, j)}{j < (2i)!(i + 1)}$ be the ascending enumeration of $e(s(i))$ and
\begin{align*}
Q_s(j, i) \dfeq & \soast{f_s(i, m)}{m \in (2i)!(j + 1) \setminus (2i)!j} \setminus \bigcup_{k, \ell < i}Q_s(\ell, k).
\intertext{Note that $\card{Q_s(j, i)} \leqslant (2i)!$ for all natural numbers $i$ and $j$. Having so defined pairwise disjoint $Q_s(j, i)$ for all $i, j < \omega$, let}
P_s(j) \dfeq & \bigcup_{i < \omega} Q_s(j, i) \text{ for all positive natural numbers } j \\
\text{and } P_s(0) \dfeq & \, \omega \setminus \bigcup_{j \in \omega \setminus 1} P_s(j), \\
\text{we have } P_s(0) \supset & \bigcup_{i < \omega} Q_s(0, i).
\end{align*}
\hfill We now consider the family of partitions $\mathcal{P} \dfeq \soast{P_s}{s \in \mathcal{E} \wedge \exists^\infty i < \omega : \card{s(i)} = (2i)!(i + 1)}$.\linebreak
As $\card{\mathcal{P}} \leqslant \card{\mathcal{E}} = \non(\meagre) < \mathfrak{s}(\mathfrak{pr})$ there is an infinite set $A$ of natural numbers which is not split by any member of $\mathcal{P}$. Let $\seq{a_i}{i < \omega}$ be the inreasing enumeration of $A$. Now consider the sequence $s \dfeq \seq{e^{-1}(\soast{a_j}{j < (2i)!(i + 1)})}{i < \omega}$. As $\mathcal{E}$ is an infinitely often equal family, there is a $t \in \mathcal{E}$ infinitely often equalling $s$. We have $P_t \in \mathcal{P}$, therefore $P_t$ does not split $A$. This means that there is a natural number $n$ such that $P_t(n) \cap A$ is finite. Let $k$ be a natural number such that $A \cap P_t(n) \subset k$ and $\ell$ be one such that $B \dfeq \soast{i \in \ell \setminus (n + 1)}{s(i) = t(i)}$ has more than $k$ elements. We have 
\begin{align*}
& A \cap \bigcup_{i \in B} Q_s(n, i) = A \cap \bigcup_{i \in B} Q_t(n, i) \subset A \cap \bigcup_{i < \ell} Q_t(n, i) \\
\subset \, & A \cap \bigcup_{i < \omega} Q_t(n, i) = A \cap P_t(n) \subset k
\end{align*}
so by the pigeonhole principle there has to be an $i \in B$ such that $Q_s(n, i) \subset \omega \setminus A$. Remember that $a_j = f_s(i, j)$ for all $j \leqslant i < \omega$. So we have
\begin{align*}
\soast{a_j}{j \in (2i)!(n + 1) \setminus (2i)!n} \setminus & \bigcup_{m, j < i} Q_s(m, j) = Q_s(n, i) \subset \omega \setminus A \text{ and hence} \\
\soast{a_j}{j \in (2i)!(n + 1) \setminus (2i)!n} \subset & \bigcup_{m, j < i} Q_s(m, j) \text{, implying } i \geqslant 1 \text{ and} \\
(2i)! \leqslant & \sum_{j, m < i} \card{Q_s(m, j)} \leqslant i \sum_{j < i} (2j)! \leqslant i^2(2i - 2)!. \\
\text{Therefore } (2i)(2i - 1) \leqslant & \, i^2, \\
\text{so } i(3i - 2) \leqslant & \, 0 \text{ contradicting both factors of this product} \\
& \quad \text{being no smaller than } 1.
\end{align*}
\end{proof}
\subsubsection{Pair-Splitting}

\begin{theo}
\label{e_ubd <= s_pair}
$\eubd \leqslant \spair$.
\end{theo}

The following proof uses an interval partition with intervals of sufficiently quickly growing lengths such that, regardless of the behaviour of a set $S$ of integers on previous intervals, a predictor will still narrow down the possibilities of that set on the next interval sufficiently to point to two elements in the current interval which are either both in $S$ or both outside of $S$.

\begin{proof}
Assume towards a contradiction that $\spair < \eubd$ and let $\mathcal{S}$ be a pair-splitting family of cardinality $\spair$. We inductively define a function $f : \omega \longrightarrow \omega$ by $f(0) \dfeq 3$ and $f(i + 1) \dfeq f(i) + 1 + 2^{2^{f(i)}}$ and use it to define an interval partition $\seq{J_i}{i < \omega}$ by setting $\displaystyle J_i \dfeq f(i) \setminus \bigcup_{k < i} J_k$ for all natural numbers $i$. Note the all intervals in this partition have odd length. Let $\seq{F_i}{i < \omega}$ be an enumeration of all finite sets of natural numbers. We define a family
\begin{align*}
\mathcal{E} & \dfeq \soast{g_S}{S \in \mathcal{S}} \text{, where} \\
g_S & \dfeq \seq{m^S_k}{k < \omega} \text{ for } S \in \mathcal{S} \text{, and} \\
m^S_k & \text{ is such that } F_{m^S_k} = S \cap J_k.
\end{align*}
As $\card{\mathcal{E}} \leqslant \card{\mathcal{S}} = \spair < \eubd$, there is a predictor $\pi = \opair{B}{\seq{\pi_i}{i \in B}}$ predicting every member of $\mathcal{E}$.
\begin{claim*}
For every natural number $\ell$ there is a pair $a_\ell = \{b_{2\ell}, b_{2\ell + 1}\} \subset J_\ell$ (where $b_{2\ell} < b_{2\ell + 1}$) such that for every $x \in \mathcal{E}$, we have $a_\ell \cap F_{\pi_\ell(x \upharpoonright \ell)} \in \{0, a_\ell\}$.
\end{claim*}
\begin{proofofclaim}
For $\ell = 0$ the claim holds true by the pigeonhole principle. In the following we show the claim to hold for $\ell = k + 1$. For any natural number $k$ there are $2^{f(k)}$ subsets of $f(k)$ and hence just as many initial segments $x \upharpoonright k$ for $x \in \mathcal{E}$. Let $\seq{t_i}{i < 2^{f(k)}}$ be an enumeration of these initial segments. Set $K_{2^{f(k)}} \dfeq J_{k + 1}$ and for every $i \leqslant 2^{f(k)}$, let $K_i$ be the bigger one of the two sets $K_{i + 1} \cap F_{\pi_k(t_i)}$ and $K_{i + 1} \setminus F_{\pi_k(t_i)}$. Inductively we have $\card{K_i} \geqslant 2^i + 1$ for every $i \leqslant 2^{f(k)}$. Now we let $b_{2k} \dfeq \min(K_0)$ and $b_{2k + 1} \dfeq \max(K_0)$. Then $a_k \dfeq \{b_{2k}, b_{2k + 1}\}$ provides what was demanded.
\end{proofofclaim}

Now we let $A \dfeq \soast{a_\ell}{\ell \in B}$. The set $A$ is an unbounded set of pairs, therefore there is an $S \in \mathcal{S}$ pair-splitting it. Let $m$ be a natural number such that $\pi_i(g_S \upharpoonright i) = g_S(i)$ for all $i \in B \setminus m$ and let $n \in B \setminus m$ be such that $a_n \in A \setminus [g(m)]^2$ is split by $S$, i.e. exactly one of $b_{2n}, b_{2n + 1}$ is in $S$, in other words, $S \cap a_n \notin \{0, a_n\}$. Then,
\begin{align*}
& \{0, a_n\} \ni a_n \cap F_{\pi_n(g_S \upharpoonright n)} = a_n \cap F_{g_S(n)} = a_n \cap F_{m^S_n} \\
= \, & a_n \cap J_n \cap S = a_n \cap S \notin \{0, a_n\} \text{, a contradiction.}
\end{align*}
\end{proof}

\begin{theo}
\label{s_pair <= ß}
$\spair \leqslant \sss$.
\end{theo}

The following proof depends on the observation that for any unbounded set $A$ of pairs of natural numbers one can define an alternating harmonic series padded with zeros (a technique employed in \cite{020BBBHHL0} and \cite{019BBH0}) such that any index set rendering the respective subseries divergent will also pair-split $A$.

\begin{proof}
Assume towards a contradiction that $\sss < \spair$ and let $\mathcal{I}$ be an $\sss$-sized family of infinite index sets of natural numbers such that for every conditionally convergent series $\displaystyle\sum_{i \nearrow \infty} r_i$ there is an index set $I \in \mathcal{I}$ such that $\displaystyle \sum_{i \in I} r_i$ is diverging. As $\card{\mathcal{I}} = \sss < \spair$, there is an unbounded collection $A$ of pairs of natural numbers which is not pair-split by any member of the family $\mathcal{I}$. Let $\seq{a_i}{i < \omega}$ be an enumeration of $A$.
\begin{align*}
B_0 \dfeq & 0; \\
B_{i + 1} \dfeq & \begin{cases}
B \cup \{a_i\} \text{ if } \max(b) < \min(a_i) \text{ for all } b \in B_i, \\
B_{i + 1} \dfeq B_i \text{ otherwise.}
\end{cases}
\end{align*}
and finally $B \dfeq \bigcup_{i < \omega} B_i$. Clearly, $B$ is an infinite unbounded subcollection of $A$, moreover we have $a_1 < a_2$ or $a_3 < a_0$ for $\big\{\{a_0, a_1\}_<, \{a_2, a_3\}_<\big\} \in [B]^2$. Let $\seq{b_i}{i < \omega}$ be the increasing enumeration of $\displaystyle \bigcup B$. \hfill Clearly, $\displaystyle B = \bsoast{\{b_{2i}, b_{2i + 1}\}}{i < \omega}$ and, being a subset of $A$, the family $B$ is not pair-split by any member of $\mathcal{I}$. Let $n(i) \dfeq \min(\soast{j < \omega}{b_j \geqslant i})$ for all natural numbers $i$. We now consider the series
\begin{align*}
\sum_{i \nearrow \infty} r_i, & \\
\text{where } r_i & \dfeq \begin{cases}
\frac{(-1)^{n(i)}}{n(i) + 1} & \text{if } i \in \bigcup B, \\
0 & \text{else.}
\end{cases}
\end{align*}
This is just the alternating harmonic series padded with zeroes, \cf\ \autoref{padding}, and therefore it is conditionally convergent. Now let $I \in \mathcal{I}$ be such that $\displaystyle \sum_{i \in I} r_i$ diverges. Let $k$ be a natural number such that $B \subset [k]^2 \cup [\omega \setminus k]^2$ and $b \subset I$ or $b \subset \omega \setminus I$ for all $b \in B \setminus [k]^2$. As $\displaystyle \sum_{i \in I} r_i$ diverges, so does $\displaystyle \sum_{i \in I \setminus k} r_i$. But $\displaystyle \sum_{i \in I \setminus k} r_i = \sum_{i \in I \cap \bigcup B \setminus k} r_i$ which, being divergent, by the Leibniz criterion cannot be an alternating series. Then $b \not \subset I$ and $b \not \subset \omega \setminus I$ for some $b \in B \setminus [k]^2$, a contradiction.
\end{proof}

Since together with Theorem \ref{e_ubd <= s_pair} we have now established the inequality $\mathfrak{e} \leqslant \sss$ between characteristics originating in, respectively, algebra and analysis, let us point out the following immediate Corollary.

\begin{coro}
For any cardinal $\kappa < 2^{\aleph_0}$ at least one of the the following holds:
\begin{itemize}
\item For all families $\mathcal{F} \in \big[[\omega]^\omega\big]^\kappa$ there is a conditionally convergent series $\displaystyle \sum_{i \nearrow \infty} r_i$ such that for all index sets $I \in \mathcal{F}$ the series $\displaystyle\sum_{i \in I} r_i$ converges.
\item There is a subgroup of the Baer-Specker group no larger than $\kappa$ exhibiting the Specker phenomenon.
\end{itemize}
\end{coro}

\FloatBarrier

\clearpage

\bibliography{t}
\bibliographystyle{t}
\end{document}